\documentclass[11pt]{article}
\usepackage{amsfonts}
\usepackage{amsmath}
\usepackage{amsbsy}
\usepackage{amssymb,latexsym}
\newtheorem{theorem}{Theorem}
\newtheorem{lemma}{Lemma}

\usepackage{graphicx}
\usepackage{epstopdf}
\newtheorem{corollary}[theorem]{Corollary}
\newtheorem{definition}{Definition}

\title{Inner tube formulas for polytopes}
\author{
 \c{S}ahin Ko\c cak\footnote{Anadolu University, Department of Mathematics,
26470  Eski\c{s}ehir, Turkey (skocak@anadolu.edu.tr)} \ and
 Andrei V. Ratiu\footnote{Istanbul Bilgi University, Department of Mathematics,
 Kurtulu\c s Deresi Cad. 47, Dolapdere, 34435 Beyoglu, Istanbul (ratiu@bilgi.edu.tr) and University of Melbourne, Department of Mathematics and Statistics, Parkville, Melbourne, VIC 3010, Australia (aratiu@unimelb.edu.au)}
}

\begin{document}
\maketitle

\begin{abstract}
We show that the
volume of the inner $r$-neighborhood of a 
polytope in the $d$-dimensional Euclidean space is a pluri-phase Steiner-like function, i.e. a continuous piecewise  polynomial function of degree $d$, proving thus a conjecture of Lapidus and Pearse. In the case when the polytope is dimension-wise equiangular we determine the coefficients of the initial polynomial as functions of the dihedral angles and the skeletal volumes of the polytope. 
We discuss also the degree of differentiability of this function
and give a lower bound in terms of the set of normal vectors of
the hyperplanes defining the polytope. We give also sufficient
conditions for the highest differentiability degree to be
attained. 
\end{abstract}

\section{Introduction}

\begin{figure}[h]
\begin{center}
\includegraphics[scale=.8]{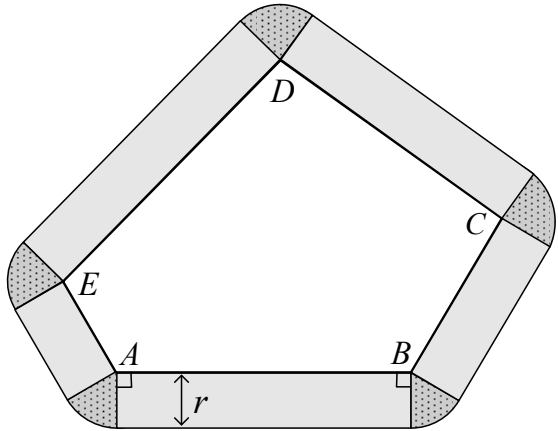}
\end{center}
\caption{The outer $r$-neighborhood of a convex polygon in the
plane. Note that the sector-volumes add up to $\pi r^2.$}
\label{outer}
\end{figure}

The Steiner formula is a beautiful cornerstone of convex
geometry dating from 1840. It states that the volume of the
$r$-parallel set of a given convex and compact set in ${\mathbb
R}^d$ can be expressed as a polynomial function of $r$ of degree
$d$ (for $r \geq 0$) (\cite{Berger} Proposition 12.3.6, \cite{Schneider} Formula 4.1.1). The $r$-parallel set $A_r$ of a set
$A\subset {\mathbb R}^d $ is by definition the Minkowski sum of
$A$ with the ball of radius $r$ in ${\mathbb R}^d$, or, in other
words, the set of points in ${\mathbb R}^d$ with distance at most
$r$ to $A$. As a simple example,
 the following formula holds for a convex, compact
set $A$ in the plane, having non-empty interior:
$$
{\rm area} (A_r) = {\rm area} (A) + {\rm length} (\partial A)\ r + \pi r^2 \qquad
(A\subset {\mathbb R}^2).
$$
The Figure \ref{outer} gives a "proof without words" for a convex polygon in the
plane.
In space, we have the following formula for a convex body
$A$, with smooth boundary $\partial A$:
$$
{\rm vol} (A_r) = {\rm vol} (A) +{\rm area} (\partial A)\ r + H(\partial A)\  r^2 +
\frac{4}{3} \pi r^3,
$$
where $H(\partial A)$ denotes the total mean curvature of
$(\partial A)$ (\cite{Morvan} Theorem 51, \cite{Gray} Theorem 10.2).

In 1939, H. Weyl proved that volumes of tubular
$r$-neighborhoods of smooth submanifolds of ${\mathbb R}^d$ can be
expressed as polynomials of $r$ and gave impressive explicit
formulas for the coefficients of the polynomial in terms of the
curvature tensor of the submanifold. As simple but interesting
examples we note two cases: For a smoothly embedded circle C
 in ${\mathbb R}^d$ it holds:
$$
{\rm vol} (C_r) = \alpha_{d-1} \ {\rm length} (C)\ r^{d-1},
$$
where $\alpha_{d-1}$ denotes the volume of the $(d-1)$-dimensional unit
ball.

For a closed orientable surface $S$ in
${\mathbb R}^d$ one has the following formula:
$$
{\rm vol} (S_r) = \alpha_{d-2} \ {\rm area} (S)\  r^{d-2} + \frac{1}{d} \
\beta_{d-1}\ \chi (S)\ r^d,
$$
where  $\beta_{d-1}$ denotes the volume of the $(d-1)$-dimensional unit
sphere and $\chi (S)$ is the Euler characteristic of $S$ (\cite{Berger} 12.10.9.2).

Now, convex sets needn't be smooth and submanifolds needn't
be convex. In 1959 Federer created the notion of "sets of positive
reach" including both of these important classes of sets, the
convex sets and (${\cal C}^2$-smooth) submanifolds in ${\mathbb R}^d$ of
any codimension. A closed set $A\subset {\mathbb R}^d$ is called
of positive reach, if there exists a parallel set $A_r$ (with
$r>0$) such that any point of $A_r$ has a unique nearest point in
$A$. The supremum of these $r$ is called the reach of $A$. Federer
proved that Steiner formula still holds for this larger class of
sets, with $r$ within the reach of the set  (\cite{Federer} Theorem 5.6).

It seems that during the long history of this problem, the inner
neighborhoods of convex bodies (or of closed hypersurfaces) was
not an issue of attraction. From the above formulas one sees that
the area of the inner $r$-neighborhood of a convex set $A$ with
smooth boundary in ${\mathbb R}^2$ can be expressed (for small
$r$) as ${\rm length} (\partial A)\ r - \pi \ r^2$ and the volume of the inner
$r$-neighborhood of a convex body with smooth boundary in ${\mathbb
R}^3$ as ${\rm area} (\partial A)\ r -  H(\partial A)\  r^2 +
\frac{4}{3} \pi r^3$.

If the boundary of a domain is not ${\cal C}^2$-differentiable or if it
is not of positive reach from inside, then there are, to our
knowledge, no available Steiner-like formulas. In the negative,
for many simple convex domains (for example the semidisc in the
plane) the volume of the inner $r$-neighborhood is definitely
non-polynomial even for small enough $r$.

The interest in volumes of inner neighborhoods of domains in
${\mathbb R}^d$ was actualized by the recent research of Lapidus
and coworkers, notably by Lapidus-Pearse, who established a
startling formula for volumes of neighborhoods of fractals,
whereby they defined new complex dimensions for fractals, related
the neighborhoods of fractals to inner neighborhoods of some
associated sequences of open domains and expressed the volume of
the $r-$neighborhood of a fractal as a sum of residues of a
certain associated $\zeta$-function at the complex dimensions of
the fractal (\cite{LapidusF} Theorem 8.1, \cite{LapidusP} Theorem 7.4, \cite{Deniz} Theorem 1).

\begin{figure}[h]
\begin{center}
\includegraphics[scale=.8]{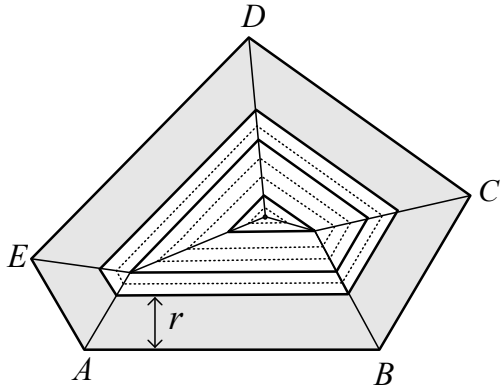}
\end{center}
\caption{The inner $r$-neighborhood of a convex polygon $P$ in the
plane. For small $r$, until the passage to a
quadrangle, its area is given by the formula $A(r)= {\rm perimeter} (P)
\cdot r - (\sum \tan (\alpha_i /2))\cdot r^2$, where the
$\alpha_i$'s are the outer angles of the polygon.} 
\label{inner}
\end{figure}

In this article we consider a simple but important type of
convex bodies, the polytopes in ${\mathbb R}^d$. The inner
neighborhoods stabilize at the inradius. We prove a conjecture of Lapidus and Pearse (\cite{LapidusP1}, Conjecture 1) stating that the volume of the inner $r$-neighborhood of a (convex, compact)
polytope is a continuous piecewise polynomial function for
$r\geq 0$.  We call such a behavior pluri-phase Steiner-like. The Figure \ref {inner}
shows a quatre-phase Steiner-like example in the plane.

We discuss also the degree of differentiability of this function
and give a lower bound in terms of the set of normal vectors of
the hyperplanes defining the polytope. We give also sufficient
conditions for the highest differentiability degree $d-1$ to be
attained. There are interesting consequences resulting from these
considerations about the geometric $\zeta$-function for the
polytopes in the sense of Lapidus-Pearse, but we want to restrict
ourselves in this article to the convex geometry framework.

\section{The volume function}
Let $P$ be a convex body in $E^d$, i.e. a compact convex subset with non-empty interior. For all positive $r$ we denote by $P(r)$ the {\em $r$-interior of} $P$:
$$P(r)=\{Q\in P | {\rm d}(Q,\partial P)\geq r \}
$$ 
and we call the set $P\setminus P(r)$ the {\em inner $r$-neighborhood of} $\partial P$. The aim of this article is the study of the inner neighborhood volume function
$$V_P(r)={\rm vol}_d(P\setminus P(r))
$$
or, equivalently, of the volume function
$$W_P(r)={\rm vol}_d(P(r))={\rm vol}_d(P)-V_P(r).$$

We start with a straightforward example.

{\bf Example} If $P$ is a $d$-dimensional Euclidean ball of radius $g$, then
$$V_P(r)=\left\{\begin{array}{lc}\sum_{i=0}^{d-1}\kappa_i(P)r^{d-i}, & \mbox{if $r\leq g$}\\  & \\{\rm vol}_d(P), & \mbox{if $r\geq g$}
\end{array}
\right.$$ Here the coefficients are
$$\kappa_i(B)=(-1)^{d-i-1} \binom{d}{i} \frac{{\rm vol}_d(B)}{g^{d-i}}, \mbox{for all $i$}.
$$

The {\em inradius} of a convex body $P$ is defined as the maximum positive value $g$ such that $P$ contains a ball of radius $g$. It is then clear that the volume function $W_P$ stabilizes at $g$, i.e.
$W_P(r)=0 \mbox{, for all $r\geq g$}
$. 

Let $d$ be a non-negative integer.
\begin{definition} {\rm A {\em degree $d$ pluri-phase Steiner-like function} $\varphi$ is a continuous function
$$\varphi : [0, +\infty)\longrightarrow [0, +\infty)
$$
such that there is a partition of the non-negative half axis into $m$ intervals $0=g_0<g_1< \dots < g_{m-1}<+\infty$ such that
\begin{itemize}
\item for each $1\leq i\leq m-1$, the restriction of $\varphi$ on the subinterval $[g_{i-1},g_i]$ is a degree $d$ polynomial
$$\varphi(r)=\sum_{k=0}^d \kappa_{d-k,i}r^k, \mbox{for all $r\in [g_{i-1},g_i]$}
$$ and
\item $\varphi$ is constant on $[g_{m-1},+\infty)$
\end{itemize}
The {\em $g$-value} of $\varphi$ is by definition $g_{m-1}$.
In case $m=2$,  $\varphi$ is called a {\em diphase Steiner-like function} and the coefficient $\kappa_{d-k,1}$ of $r^k$ is denoted simply by $\kappa_{d-k}$.} 
\end{definition}

The main result of this article is to show that whenever $P$ is a polytope in $E^d$, i.e. a convex body that is the intersection of finitely many closed half-spaces in $E^d$, the inner neighborhood volume function $V_P$ is a degree $d$ pluri-phase Steiner-like function. If $\{{\cal H}_1, {\cal H}_2, \dots , {\cal H}_m\}$ is a minimal collection of half-spaces in $E^d$ defining $P$, the sets $F_j= P \cap \partial {\cal H}_j$ are called the {\em facets} of $P$. Their union equals the boundary of the polytope. We will use the following differentiation formula for the volume function:
$$W'_P(0)=-\sum_{j=1}^m{\rm vol}_{d-1}(F_j),
$$ valid for all $d$-dimensional polytope $P$ with facets $F_1, F_2, \dots , F_m$. As each facet $F_j$ is a $(d-1)$-dimensional polytope, the proof will essentially consist of an inductive argument on dimension. By definition, a $0$-dimensional polytope is a point and its $0$-dimensional volume is simply $1$.

By contrast, even for very simple convex subsets in the plane the inner neighborhood volume function is not necessarily a Steiner-like function, as the following example shows.

{\bf Example} Let $P=\{(x,y)\in R^2 | 0\leq y \leq \sqrt{1-x^2}\}$. Then the $r$-interior of $P$ is
$$P(r)=\{(x,y)\in R^2 | r\leq y \leq \sqrt{(1-r)^2-x^2}\}
$$
and 
$$V_P(r)=\left\{\begin{array}{lc}\frac{\pi}{2}-(1-r)^2\arccos \left( \frac{r}{1-r}\right)+r\sqrt{1-2r}, & \mbox{if $r\leq \frac{1}{2}$}\\  & \\ \frac{\pi}{2}, & \mbox{if $r\geq \frac{1}{2}$}
\end{array}
\right.$$ 
Note that $V_P$ is of differentiability class ${\cal C}^1$ on $[0, +\infty)$, but not ${\cal C}^2$, as:
$$\lim_{r \nearrow \frac{1}{2}}V_P''(r)=-\infty
$$

The differentiability class of the volume function reflects the metric properties of the polytope. For instance among all the rectangles, only the cubes have volume function of maximum differentiability class, as the following example proves it:

{\bf Example} Let $a_1 \leq a_2 \leq \dots \leq a_d$ and let $R$ be the $d$-dimensional rectangle 
$$R=R_{a_1,a_2,\dots ,a_d }=\{x=(x_1, x_2, \dots , x_d) | \left\vert x_i \right\vert \leq a_i, \mbox{for all $i$}\}
$$ of inradius $a_1$. Then $V_R$ is a diphase Steiner-like function:
$$V_R(r)=\left\{\begin{array}{lc}
 2^da_1a_2\dots a_d-2^d(a_1-r)(a_2-r)\dots (a_d-r),
& \mbox{if $r\leq g$}\\ & \\ 2^da_1a_2\dots a_d, & \mbox{if $r\geq g$}
\end{array}
\right.$$ 
Thus the coefficients are 
$$\kappa_{d-k}=(-1)^{k-1}{\rm vol}_{d-k}(R_{(d-k)}), \mbox{for $1\leq k \leq d$}.
$$ Here $R_{(d-k)}$ denotes the $(d-k)$-dimensional skeleton of the polytope $R$, i.e. the union of all the $(d-k)$-dimensional faces of $R$. In particular, $\kappa_{0}=(-1)^{d-1}\sharp(R_{(0)})
$. Note also the important fact that, if $m= \sharp \{1\leq i \leq d|a_i=g\}$, then $V_R$ is of differentiability class ${\cal C}^{(m-1)}$, but not ${\cal C}^m$. Indeed, for all $r$ such that $0\leq r \leq g$, $g$ is a root of the polynomial $2^d(a_1-r)(a_2-r)\dots (a_d-r)$ of multiplicity exactly $m$. 

The following theorem characterizes the polytopes in $E^d$ for which the function $V_P$ is a diphase Steiner-like function of maximal differentiability class, i.e. of class ${\cal C}^{(d-1)}$.

\begin{theorem}
A necessary and sufficient condition for the inner neighborhood volume function $V_P$ of a $d$-dimensional polytope $P$ to be a diphase Steiner-like function of class ${\cal C}^{(d-1)}$ on $[0,+\infty)$ is that $P$ admits an inscribed Euclidean $d$-dimensional ball. Moreover, in this case the coefficients $\kappa_i(P)$ of $P$ satisfy:
$$\kappa_i(P)=(-1)^{d-i-1}\binom{d}{i} \frac{{\rm vol}_d(P)}{g^{d-i}} ,
$$
for all $i=0,1,\dots, d-1$.
\label{inscribed}

\end{theorem}

{\bf Proof} We first prove the sufficiency of the condition. Assume that the origin $O$ is the center of a $d$-dimensional ball of radius $g$ inscribed in $P$, i.e. tangent to all the facets of $P$. For each $r$ in $[0,g]$ let $\tilde{P}(r)=(1-\frac{r}{g})P$. Note that, due to the assumption above, the polytopes $\tilde{P}(r)$ and $P(r)$ coincide, for all $r$. Thus for each $r$ in $[0,g]$:
$$V_P(r)={\rm vol}_d(P\setminus \tilde{P}(r))={\rm vol}_d(P)\left(1-\left(1-\frac{r}{g}\right)^d\right)={\rm vol}_d(P)\sum_{i=0}^{d-1}(-1)^{d-i-1} \binom{d}{i} \left(\frac{r}{g}\right)^{d-i} 
$$
and for $r>g$, $V_P(r)={\rm vol}_d(P)$.

In order to prove the necessity, note that for any $d$-dimensional polytope $P$, the inradius $g$ satisfies
$$g=\min\{r>0|{\rm vol}_d(P(r))=0\}=\min\{r>0|V_P(r)={\rm vol}_d(P)\}.
$$ Assume now that $V_P$ is diphase Steiner-like function. Then its $g$-value equals the inradius of $P$. Writing down the continuity conditions at $r=g$ for the function 
$$V_P(r)=\left\{\begin{array}{cc}\sum_{i=0}^{d-1}\kappa_i(P)r^{d-i} & \mbox{if $r<g$}\\ & \\ {\rm vol}_d(P) & \mbox{if $r\geq g$}
\end{array}
\right.$$ and for its first $d-1$ derivatives we obtain the triangular homogeneous system of linear equations in the coefficients $\kappa_i$:
$$\sum_{i=0}^{d-j}\binom{d-i}{j}\kappa_i(P)g^{d-i-j}=0,  \mbox{for all $j=0,\dots, d-1$}
$$
whose only family of solutions is the one stated in the theorem. Here, for simplicity we have denoted ${\rm vol}_d(P)$ by $-\kappa_d(P)$.

Consider that the ball of radius $g$ centered at the origin is contained in $P$. Then $(1-\frac{r}{g})P\subset P(r)$. Since 
$$\begin{array}{rl}
 V_P(r) & =\sum_{i=0}^{d-1}\kappa_i(P)r^{d-i}\\ & \\
      &={\rm vol}_d(P)\sum_{i=0}^{d-1}(-1)^{d-i-1} \binom{d}{i} \left(\frac{r}{g}\right)^{d-i}  \\ & \\
      & ={\rm vol}_d(P)\left(1-\left(1-\frac{r}{g}\right)^d\right)=
{\rm vol}_d\left(P\setminus \left(1-\frac{r}{g}\right)P\right),
\end{array}
$$
 we conclude that $P\setminus P(r)=P\setminus (1-\frac{r}{g})P$, for all $r<g$. 
Thus the distance of the origin to all the facets is exactly $g$, and the ball of radius $g$ considered is inscribed in $P$. qed.

\section{The equiangular polytope case}

Considering that the volume function of a convex body is a polynomial for at least a small interval $[0, \varepsilon]$, one is interested in the geometric information carried by the coefficients of this polynomial. This section answers this question explicitly for the special case of dimension-wise equiangular polytopes. 

\begin{definition} {\rm  A polytope $P$ is {\em equiangular} if all the dihedral angles between pairs of adjacent facets of $P$ are equal. For an equiangular $d$-dimensional polytope denote by $\alpha_d$ its {\em outer dihedral angle value}, i.e. for any two different facets $F_i$ and $F_j$, meeting along a common ridge, and having inner unit normal vectors $N_i$ and $N_j$, $\langle N_i, N_j \rangle= \cos \alpha_d$ holds.} \end{definition}

\begin{definition} {\rm  A polytope $P$ is {\em dimension-wise equiangular} if there exists constants $\alpha_2, \alpha_3, \dots , \alpha_{d}\in (0, \pi)$ such that any $k$-dimensional face $F$ of $P$ is an equiangular polytope with outer dihedral angle $\alpha_{k}$, for all $k$ such that  $2 \leq k \leq d$.} 
\end{definition}

As a consequence of the definitions above, all $1$-dimensional polytopes are dimension-wise equiangular.
 
{\bf Example} Following \cite{Berger} (Definition 12.5.1), a {\em flag} of a polytope $P$ is a $d$-tuple $(F_0, F_1, \dots , F_{d-1})$, consisting of $i$-faces of $P$ and $F_i \subset F_{i+1}$ for all $i$ such that $0\leq i \leq d-2$. A polytope $P$ is called {\em regular} if the symmetry group of $P$, i.e. the group of all Euclidean isometries mapping $P$ onto itself, acts transitively on the flags of $P$. As it is clear that a regular polytope $P$ in $E^d$ is equiangular and as every $k$-face of $P$ is regular, we conclude that every regular polytope is dimension-wise equiangular. Moreover, any of its $k$-faces has same inradius $\gamma_k g$, where $g$ is the inradius of $P$, $\gamma_d$ is obviously $1$, and 
$$\gamma_k=\tan\left(\frac{\alpha_{k+1}}{2}\right) \tan \left(\frac{\alpha_{k+2}}{2}\right) \dots \tan \left(\frac{\alpha_{d}}{2}\right), \mbox{for all $k$ such that $1\leq k \leq d-1$}.
$$ The expressions $\gamma_k$, defined as functions of the outer dihedral angles, have a geometric meaning even in the case of a non-regular but dimension-wise equiangular polytope, as we will see in the theorem \ref{equiangular}.

\begin{definition} {\rm  Let $P$ be a $d$-dimensional polytope with facets $F_1, F_2, \dots , F_m$. For $0\leq k \leq d$, we define inductively:
$$\Omega_k(P)=\left\{\begin{array}{lc}{\rm vol}_d(P), & \mbox{if $k=d$}\\ & \\\sum_{i=1}^m\Omega_k(F_i), & \mbox{otherwise}
\end{array}
\right.
$$. }
\label{Omega}
\end{definition}
 
 Let $F'_1, F'_2, \dots , F'_n$ be all the $k$-faces of $P$, and suppose that, for all $i$, the polytope $P$ admits exactly $\mu_i$ flags on $F'_i$, i.e. $(d-k)$-tuples of $i$-faces $(F''_k, F''_{k+1}, \dots , F''_{d-1})$, with $F''_k=F'_k$ and $F''_i\subset F''_{i+1}$, for all $i$ such that $k\leq i \leq d-2$. By a straightforward induction on $d$, one can show that 
 $$\Omega_k(P)=\sum_{i=1}^n\mu_i {\rm vol}_k(F'_i)
  .$$

\begin{theorem}
\label{equiangular}
Let $P$ be a dimension-wise equiangular polytope in $E^d$, with outer dihedral angles $\alpha_2, \dots , \alpha_{d}$. Then there exists a positive $\varepsilon$ such that, on the interval $[0,\varepsilon]$, the volume function $W_P(r)={\rm vol}_d(P(r))$ is a polynomial function and 
$$W_P(r)=\sum_{k=0}^d(-1)^{d-k}\Omega_{k}(P)\gamma_{k+1}\gamma_{k+2}\dots \gamma_{d}\frac{r^{d-k}}{(d-k)!},
$$ for all $r$ such that $0\leq r \leq \varepsilon$, where $\gamma_k$ are given by the formulas above.
\end{theorem}

{\bf Proof} By induction on the dimension $d$. Denote by $g$ the inradius of $P$. For $d=1$, take $\varepsilon=\frac{g}{2}$. 

Assume the theorem has been proved for $d-1$. Denote by $F_1, F_2, \dots , F_m$ the facets of a dimension-wise equiangular polytope $P$ in $E^d$. By the induction hypothesis, there exists an $\varepsilon' $ such that each volume function $W_{F_i}$ is a degree $d-1$ polynomial function on $[0,\varepsilon']$ and satisfies 
$$W_{F_i}(r)=\sum_{k=0}^{d-1}(-1)^{d-k-1}\Omega_{k}(F_i)\gamma'_{k+1}\gamma'_{k+2}\dots \gamma'_{d-1}\frac{r^{d-k-1}}{(d-k-1)!},
$$ for all $r$ such that $0\leq r \leq \varepsilon'$. Here 
$$\gamma'_{l}=\tan\left(\frac{\alpha_{l+1}}{2}\right) \tan \left(\frac{\alpha_{l+2}}{2}\right) \dots \tan \left(\frac{\alpha_{d-1}}{2}\right) \mbox{, for all $l$ such that $1\leq l \leq d-1$.}
$$
Take 
$$\varepsilon=\min\left\{g, \frac{\varepsilon'}{\tan\left(\frac{\alpha_{d}}{2}\right) }\right\}.
$$

For all $r$ in the interval $[0,\varepsilon]$, $P(r)$ is a dimension-wise equiangular polytope with facets $F_1(\tan\left(\frac{\alpha_{d}}{2}\right)r),$ $F_2(\tan\left(\frac{\alpha_{d}}{2}\right)r), \dots , F_m(\tan\left(\frac{\alpha_{d}}{2}\right)r)$. Applying the differentiation formula to the volume function  for the polytope $P(r)$, we get:
$$\begin{array}{rl}
 W'_P(r)& =-\sum_{i=1}^m W_{F_i}\left(\tan\left(\frac{\alpha_{d}}{2}\right)r\right) \\ & \\
      & =\sum_{i=1}^m\sum_{k=0}^{d-1}(-1)^{d-k}\Omega_{k}(F_i)\gamma'_{k+1}\gamma'_{k+2}\dots \gamma'_{d-1}\left(\tan\left(\frac{\alpha_{d}}{2}\right)\right)^{d-k-1}\frac{r^{d-k-1}}{(d-k-1)!}\\ & \\
      & =\sum_{k=0}^{d-1}\sum_{i=1}^m(-1)^{d-k}\Omega_{k}(F_i)\gamma_{k+1}\gamma_{k+2}\dots \gamma_{d}\frac{r^{d-k-1}}{(d-k-1)!}\\ & \\
       & =\sum_{k=0}^{d-1}(-1)^{d-k}\Omega_{k}(P)\gamma_{k+1}\gamma_{k+2}\dots \gamma_{d}\frac{r^{d-k-1}}{(d-k-1)!}.\\
\end{array}
$$
The theorem follows, since $W_P(0)={\rm vol}_d(P)=\Omega_d(P).$ qed

As a consequence of Theorems \ref{equiangular} and \ref{inscribed} and of the equality stated after the definition \ref{Omega}, we get the following corollary,  where $P_{(k)}$ denotes the $k$-dimensional skeleton of $P$, i.e. the union of all its $k$-dimensional faces.
\begin{corollary} 
\label{equiangular2} 
Let $P$ be a dimension-wise equiangular polytope in $E^d$ such that for every $k\geq 1$, the number of  $k$-faces of $P$ meeting at any $(k-1)$-face of $P$ is exactly $\mu_{(k)}$. Then there exists a positive $\varepsilon$ such that, on the interval $[0,\varepsilon]$, the volume function $W_P(r)={\rm vol}_d(P(r))$ is a polynomial function and 
$$W_P(r)=\sum_{k=0}^d(-1)^{d-k}\mu_{(k+1)}\mu_{(k+2)}\dots \mu_{(d)}{\rm vol}_{k}(P_{(k)})\gamma_{k+1}\gamma_{k+2}\dots \gamma_{d}\frac{r^{d-k}}{(d-k)!},
$$ for all $r$ such that $0\leq r \leq \varepsilon$.
In particular, if $P$ is a $d$-dimensional regular polytope, for all $k$ such that $0\leq k \leq d-1$,
$$d!{\rm vol}_{d}(P_{(d)})=k!\mu_{(k+1)}\mu_{(k+2)}\dots \mu_{(d)}{\rm vol}_{k}(P_{(k)})\gamma_{k+1}\gamma_{k+2}\dots \gamma_{d}g^{d-k} .
$$
\end{corollary}

{\bf Example} Consider a regular dodecahedron of unit edge length. Cut it along two planes, both parallel to a couple of opposite sides, each one being at distance $\delta=\sqrt{1-\frac{2}{\sqrt{5}}}$ from the corresponding side. We obtain a dodecahedron $D$ with all outer dihedral angles equal to $\arctan 2$, and whose two sides are regular pentagons of side $\frac{2}{\varphi}$, and remaining 10 sides are equiangular pentagons of edges 1, 1, $\frac{1}{\varphi}$, $\frac{1}{\varphi}$ and $\frac{2}{\varphi}$. Here $\varphi=\frac{\sqrt{5}+1}{2}$ is the golden ratio.  Thus the area of each of the small sides is $\frac{\sqrt{15-5\varphi}}{2}$ and of each of the 2 larger sides is $\sqrt{15-5\varphi}$. The constants are
$$\begin{array}{rl}
\gamma_3 & =1\\ & \\
\gamma_2 & =\varphi -1\\ & \\
\gamma_1 & =\sqrt{18-11\varphi}.\\
\end{array}$$

We obtain

$$\begin{array}{rl}
 W_D(r)& =\sum_{k=0}^3(-1)^{3-k}\mu_{(k+1)}\mu_{(k+2)}\dots \mu_{(d)}{\rm vol}_{k}(P_{(k)})\gamma_{k+1}\gamma_{k+2}\dots \gamma_{3}\frac{r^{3-k}}{(3-k)!} \\ & \\
      & =-20\sqrt{47-29\varphi}r^3+(50-20\varphi)r^2-7\sqrt{15-5\varphi}r+{\rm vol}_3(D),\\
\end{array}
$$
for all $r$ such that $0\leq r \leq \varepsilon$.

\begin{figure}[h]
\begin{center}
\includegraphics[scale=1]{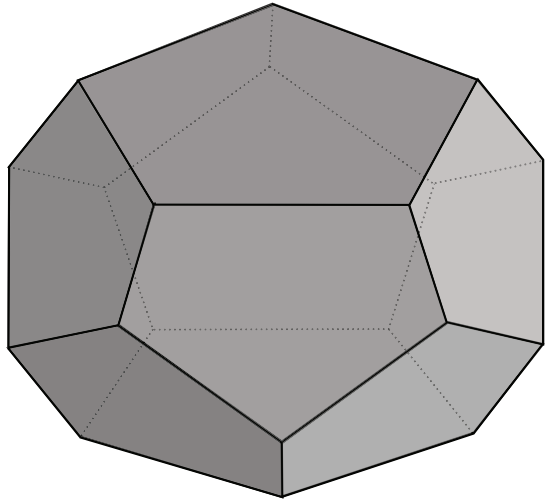}
\end{center}
\caption{ A non-regular but dimension-wise equiangular dodecahedron} \label{dodec}
\end{figure}

\section{The main result}
In this section we prove that the volume function of any polytope is a pluri-phase Steiner-like function. For that we introduce the concept of a gliding arrangement of hyperplanes in a Euclidean space. The hyperplanes supporting the facets of the "shrinking" polytope $P(r)$, as $r$ increases, constitute the gliding arrangement of hyperplanes adapted to the polytope $P$. This concept will allow an inductive argument on the dimension of the polytope.

By a hyperplane in $E^d$ we will understand an oriented hyperplane $H$, where the orientation is given by a choice of a unit normal vector $N$.
\begin{definition} {\rm 
A {\em $d$-arrangement of hyperplanes} is a finite family ${\cal A}=\{H_1,H_2,\dots ,H_m\}$ of (not necessarily distinct) hyperplanes in a $d$-dimensional real affine Euclidean space $E^d$.}
\end{definition}

\begin{definition} {\rm 
The {\em signed distance} ${\rm d}(Q,H_0)$ from a point $Q$ in $E^d$ to the oriented hyperplane $H_0$ is by definition ${\rm d}(Q,H_0)=\langle \overline{Q_0Q}, N_0 \rangle$, where $Q_0$ is any point in the hyperplane $H_0$.
The {\em cells} $C_{\varepsilon}$ of the arrangement ${\cal A}$ are indexed by $\varepsilon=(\varepsilon_1, \varepsilon_2,\dots ,\varepsilon_m)\in\{\pm 1\}^m$, and defined as the set 
$$C_{\varepsilon}=\{Q\in E^d | \varepsilon_j {\rm d}(Q,H_j) \geq 0, \mbox{ for all $j$}\}.
$$ }
\end{definition}
By convention, we consider cells with different multi-indices as different, even when empty. So a $d$-arrangement of cardinality $m$, defines exactly $2^m$ cells in $E^d$.

\begin{definition} {\rm 
We say that the family of vectors ${\cal N}=\{N_1, N_2, \dots ,N_m\}$ in $R^d$ has {\em absolute rank} $k$ if:
(i) any $k$ elements of ${\cal N}$ are linearly independent, and
(ii) there are $k+1$ vectors in the family that are linearly dependent.
The {\em rank} (respectively the {\em absolute rank}) of a $d$-arrangements ${\cal A}$ is defined to be the rank (respectively, the absolute rank) of the family of the normal vectors ${\cal N}$ of its elements.}
\end{definition}

The chain of inequalities holds: $1 \leq k \leq {\rm rk}{\cal N} \leq \min\{m,d\}$, for all non-empty families ${\cal N}$. By convention, we consider the absolute rank of an empty family to be null.

\begin{definition} {\rm 
Let $H_0$ be an oriented hyperplane in $E^d$, with $d\geq 2$. The {\em trace} of  the $d$-arrangement ${\cal A}$ on $H_0$ is the family $\tilde{\cal A}_0$ of all the proper intersections of $H_0$ with each of  the hyperplanes in ${\cal A}$, each intersection being endowed with the normal vector $\tilde{N}_j=\frac{pr_0(N_j)}{|pr_0(N_j)|}$, where $pr_0(v)=v-\langle v,N_0\rangle N_0$ is the projection on the linear hyperplane subtending $H_0$, i.e.
$$\tilde{\cal A}_0=\{(H_j\cap H_0,\frac{pr_{0}(N_j)}{|pr_{0}(N_j)|})|1\leq j \leq m \mbox{ such that $H_j\cap H_0$ is a proper subspace in $H_0$}\}.
$$}
\end{definition}

Note that, for all $j$ such that $1\leq j\leq m$, the trace $\tilde{\cal A}_j$ on the hyperplane $H_j$ is a (possibly empty) $(d-1)$-arrangement of absolute rank at least $k-1$ and of rank at least ${\rm rk}{\cal A}-1$. In particular, if ${\rm rk}{\cal A}\geq 2 $, then all the traces $\tilde{\cal A}_j$ of ${\cal A}$ on its own hyperplanes are non-empty. 

Consider now two hyperplanes $H_a$ and $H_b$. Note that, for all $j$ such that $1\leq j\leq m$, if the intersection $H_j\cap H_a \cap H_b$ is a proper affine subspace in $H_a\cap H_b$, then the projection $pr_{ab}(N_j)$ of $N_j$ on the linear subspace subtending $H_a\cap H_b$ cannot vanish. Moreover the trace of the arrangement $\tilde{\cal A}_a$ on $H_a\cap H_b$ (as a hyperplane in $H_a$) equals
$$\{(H_j\cap H_a \cap H_b,\frac{pr_{ab}(N_j)}{|pr_{ab}(N_j)|})|1\leq j \leq m \mbox{ such that $H_j\cap H_a \cap H_b$ is a proper subspace in $H_a \cap H_b$}\}
$$
so it coincides with the trace of $\tilde{\cal A}_b$ on the same $H_a\cap H_b$ (as a hyperplane in $H_b$). We say thus that the iterated trace of an arrangement ${\cal A}$ on several hyperplanes does not depend on the particular order in which these traces are taken. 

\begin{definition} {\rm 
A {\em gliding $d$-arrangement of hyperplanes} is a family 
$${\cal A}=\{(H_1,v_1),(H_2,v_2),\dots ,(H_m,v_m)\}$$ where each hyperplane $H_j$ is endowed with a normal velocity vector $v_j$.  A pair $(H,v)$ as above will be called a {\em gliding hyperplane}. }
\end{definition}

{\bf Example} Let $P$ be a $d$-dimensional polytope in $E^d$. Consider the gliding $d$-arrangement  ${\cal A}_P$ consisting of all the hyperplanes supporting the facets of $P$ with the orientation given by the inner normal vectors $N_j$ and with velocity vectors $v_j=N_j$ for all $j$. We call ${\cal A}_P$ the gliding arrangement {\em adapted to} $P$. 

A gliding arrangement ${\cal A}$ determines at each moment of time $t$ a $d$-arrangement $${\cal A}(t)=\{H_1(t),H_2(t),\dots ,H_m(t)\},$$ where
$H_j(t)$ is the image of $H_j(0)=H_j$ under the translation by vector $tv_j$, i.e.:
$$H_j(t)=\{Q\in E^d | {\rm d}(Q,H_j)-t\langle v_j,N_j\rangle=0\}
$$
 Note that the rank and the absolute rank of ${\cal A}(t)$ are constant. For a multi-index $\varepsilon \in \{\pm 1\}^m$ denote the  
cell $C_{\varepsilon}$ by $C(t)$ to indicate its dependance on the time parameter. Precisely:
$$C(t)=\{Q\in E^d | \varepsilon_{\rm d}(Q,H_j(t))\geq 0, \mbox{for all $j=1,\dots,m$}\}
$$
Note that if $C(t_0)$ is bounded for some $t_0$ then it is bounded for all real $t$. Let 
$$W(t)=W_{\varepsilon}(t)={\rm vol}_d(C(t)).
$$ be its $d$-dimensional volume function. 
\begin{definition} {\rm  A gliding $d$-arrangement ${\cal A}$ is {\em of first type} if it admits a bounded cell $C_{\varepsilon}(t_0)$ with non-zero $d$-dimensional volume, for some multi-index $\varepsilon$ and a real $t_0$.}
\end{definition}

{\bf Example} The gliding arrangement adapted to a polytope is of first type by definition.

\begin{lemma}
Let ${\cal A}$ be a gliding arrangement of hyperplanes in $E^d$ of first type. Let $C(t)$ be a bounded cell of ${\cal A}$ of volume function $W$. Then $W$ is continuous and ${\rm supp}\,W$ is an interval. For all $t\not\in {\rm supp}\,W$, $C(t)$ is empty. 
\label{support}
\end{lemma}

{\bf Proof} The continuity of $W$ follows simply by induction on the dimension $d$ and we will skip its proof. 
If $Q_i\in C(t_i)$ for $i=0,1$, then for all $\lambda\in [0,1]$, the point $(1-\lambda)Q_0+\lambda Q_1\in C((1-\lambda)t_0+\lambda t_1)$. Thus if $t_0\neq t_1$ and if ${\rm vol}_d(C_{\varepsilon}(t_1))>0$, then ${\rm vol}_d(C((1-\lambda)t_0+\lambda t_1)) \geq \lambda^d {\rm vol}_d(C(t_1))>0$. Thus if $t_0 < t_1$ are in ${\rm supp}\,W$, then the whole interval $[t_0,t_1]$ is contained in ${\rm supp}\,W$. By the same argument as above, if $t_1\in {\rm supp}\,W$, and $t_0\not\in {\rm supp}\,W$, but if $C(t_0)$ is non-empty then $t_0$ is in the closure of ${\rm supp}\,W$, which is absurd. Thus $C(t)$ is empty for all $t \not\in 
{\rm supp}\,W$. qed.

Assume $(H_j,v_j)$ and $(H_0,v_0)$ are gliding hyperplanes in $E^d$ in general position.
\begin{definition} {\rm 
We define the {\em trace of the gliding hyperplane $H_j$ on $H_0$} to be the hyperplane $H_j \cap H_0$ endowed with the normal vector $\tilde{N}_j=\frac{pr_0(N_j)}{|pr_0(N_j)|}$ and with the velocity vector 
$$\tilde{v}_j=\langle v_j-v_0, N_j\rangle \frac{pr_0(N_j)}{|pr_0(N_j)|^2}$$.}
\end{definition}
Note that in the orthogonal complement of $H_0 \cap H_j$, $\{N_0, N_j\}$ and $\{\frac{pr_j(N_0)}{|pr_j(N_0)|^2},\frac{pr_0(N_j)}{|pr_0(N_j)|^2}\}$ are dual bases. The vector $\tilde{v}_j$ is exactly the $ \frac{pr_0(N_j)}{|pr_0(N_j)|^2}$-component of the vector $v_j-v_0$.

Let ${\cal A}=\{(H_1,v_1),(H_2,v_2),\dots ,(H_m,v_m)\}$ a gliding $d$-arrangement of hyperplanes and $(H_0,v_0)$ be a gliding hyperplane in $E^d$, with $d\geq 2$. 
\begin{definition} {\rm 
We define the {\em trace of the gliding $d$-arrangement ${\cal A}$ in $H_0$} to be the family $\tilde{\cal A}_0$ of all the traces of the gliding hyperplanes $(H_j,v_j)$ on $(H_0,v_0)$, where we choose only those hyperplanes $(H_j,v_j)$ such that $H_j$ and $H_0$ are in general position.}
\end{definition}

Let ${\cal A}_P$ be the gliding arrangement adapted to some $d$-dimensional polytope $P$ in $E^d$. Note that, if $P$ is a $d$-dimensional rectangle, the traces of ${\cal A}_P$ on the hyperplanes $H_j$ containing the facets of $P$ are adapted to the respective facet $P\cap H_j$. Nevertheless, for a general polytope $P$ this is not necessarily true.

\begin{lemma}
Let $H_0$ be a gliding hyperplane and ${\cal A}$ a gliding arrangement of hyperplane in $E^d$, with $d\geq 2$, and denote by $\tilde{\cal A}_0$ the trace of ${\cal A}$ on $H_0$. Then the arrangement $\tilde{\cal A}_0(t)$ is the trace of ${\cal A}(t)$ on $H_0(t)$, for all $t$.
\end{lemma}

In particular, we conclude that the iterated trace of a gliding arrangement ${\cal A}$ on several gliding hyperplanes does not depend on the particular order in which these traces are taken.

{\bf Proof} We can assume that ${\cal A}$ consists of a single gliding hyperplane $(H,v)$ and that the normal vectors $N$ and $N_0$ are linearly independent. We can identify $E^d$ with the Euclidean real coordinate space ${\mathbb R}^d$ in such a way that $N_0={\mathbf e}_1$, $N=\frac{1}{\sqrt{b^2+c^2}}(b{\mathbf e}_1+c{\mathbf e}_2)$, with $c\neq 0$, $H_0(0)=\{x_1=0\}$ and $H(0)=\{bx_1+cx_2=0\}$.

Then $v_0=a{\mathbf e}_1$ and $v=\omega (b{\mathbf e}_1+c{\mathbf e}_2)$, for some real $\omega$, thus $H_0(t)=\{x_1=at\}$ and $H(t)=\{bx_1+cx_2=\omega(b^2+c^2)t\}$ intersect at 
$$\{(at,\frac{\omega(b^2+c^2)-ab}{c}t,x_2, \dots , x_d) | t\in {\mathbb R}\}.
$$
Thus the trace $\widetilde{(H,v)}_0$ of the gliding hyperplane $(H,v)$ on $H_0$ has velocity vector $\frac{\omega(b^2+c^2)-ab}{c}{\mathbf e}_2=
\langle v-v_0, N\rangle \frac{pr_0(N)}{|pr_0(N)|^2}$. The Lemma is thus proven.

\begin{theorem}
Let ${\cal A}$ be a gliding $d$-arrangement of hyperplanes in $E^d$ of first type and of absolute rank $k$ and let $C$ be a bounded cell, such that for some value $t_0$, $C(t_0)$ has non-zero $d$-volume. Then the function $W(t)={\rm vol}_d(C(t))$ is a pluriphase degree $d$ Steiner-like function of class ${\cal C}^{k-1}$ on ${\mathbb R}$. 
\label{pluriphase}
\end{theorem}

{\bf Proof} By induction on the dimension $d$.

The theorem is clear if $d$ is 1. Note that in this case $W(t)$ is of class ${\cal C}^1$ only if the cell is non-empty for all $t$.

Assume the theorem has been proven for $(d-1)$-arrangements. The cell $C=C_{\varepsilon}$ of the $d$-arrangement ${\cal A}$ determines in the trace $\tilde{\cal A}_j$ a unique cell $C_j$ for each $j$ and, by the Lemma above, $C_j(t)$ is the $j$th facet of $C(t)$, for all $t$ in ${\rm supp}\,W$. Since $C$ is a bounded cell, the same is true for $C_j$. In particular we conclude that $\tilde{\cal A}_j$ is a non-empty $(d-1)$-arrangement for all $j$, as the empty arrangement has only one cell, which is unbounded. Note, however, the important fact that if ${\rm vol}_{d-1}C_j(t_1)\neq 0$ for at least one value $t_1$ for which ${\rm vol}_d(C(t_1))=0$, then ${\cal A}$ must contain at least one hyperplane $H_l$, with $l\neq j$ such that the normal vectors $N_l$ and $N_j$ are linearly dependent. We conclude that ${\cal A}$ has absolute rank $1$ in this case.

Denote the $(d-1)$-dimensional volume function of $C_j(t)$ by $W_j(t)$, for each $j$. Note also that the inner unit normal vector of the cell $C(t)$ along its facet $C_j(t)$ is $\varepsilon_j N_j$. By taking the one-sided derivatives with respect to $t$, we get
$$W'(t\pm)=-\chi (t\pm)\sum_{j=1}^m \langle v_j,\varepsilon_j N_j \rangle W_j(t\pm), \mbox{for all $t$}
$$
where $\chi $ is the characteristic function of the interval ${\rm supp}\,W$. 

Since for each $j$, the trace $\tilde{\cal A}_j$ is a non-empty $(d-1)$-arrangement of absolute rank at least $k-1$, we conclude that the function $W$ is a pluriphase degree $d$ Steiner-like function of class ${\cal C}^{k-1}$ on ${\mathbb R}\backslash \partial ({\rm supp}\,W(t))$.

If the absolute rank of ${\cal A}$ is $1$, the conclusion of the theorem is proven. Otherwise, the absolute rank of ${\cal A}$ is greater than $1$, and by the argument above, the volume function $W_j$ is null outside the interval  ${\rm supp}\,W$, for all facet $j$. The theorem follows. qed

\begin{corollary}
Let $P$ be a $d$-dimensional polytope in ${\mathbb R}^d$ and let ${\cal N}=\{N_1, N_2,..., N_m\}$ be the inner normal vectors of the facets of P.
If the family of normal vectors ${\cal N}$ has absolute rank $k$, then the inner neighborhood volume function $V_P(r)$ is a pluriphase degree $d$ Steiner-like function of class ${\cal C}^{k-1}$ on $[0,+\infty)$.
\label{convexinner}
\end{corollary}

{\bf Proof} Consider the gliding $d$-arrangement ${\cal A}_P$ adapted to $P$. The polytope $P$ is the cell $C=C_{(+1,+1,\dots,+1)}$. We first show that for non-negative $t$, $C(t)=C_{(+1,+1,\dots,+1)}(t)$ coincides with the $t$-interior of $P$, i.e.:
$$P(t)=\{Q\in P|{\rm d}(Q,\partial P)\geq t\}
$$
Indeed, if $Q\in C(t)$ then for all $j$, ${\rm d}(Q,H_j(t))\geq 0$. So ${\rm d}(Q,P\cap H_j)\geq {\rm d}(Q,H_j)\geq t$, for all facet $P\cap H_j$ of $P$. 
Conversely, if $Q\in P$ is a point satisfying ${\rm d}(Q,\partial P) \geq t$ then ${\rm d}(Q,\overline{{\mathbb R}^d \setminus P}) \geq t$. Since any hyperplane $H_j(0)$ supporting a facet of $P$ is contained in the complement $\overline{{\mathbb R}^d \setminus P}$, we conclude that ${\rm d}(Q,H_j(0)) \geq t$ for all $j$. Thus $Q\in C(t)$. The proof follows now by directly applying the theorem to the gliding $d$-arrangement ${\cal A}_P$ and to the cell $P$ and the relation $W(r)={\rm vol}_d(P)-V_P(r)$, valid for all positive $r$. qed.

Moreover, the theorem below proves that, for each $0 \leq k \leq d$, there exists a polytope $P$ such that the absolute rank of ${\cal N}$ is $k$ but $V_P(r)$ is not of class ${\cal C}^{k}$ on $[0,+\infty)$. So, in this sense, the bound for the differentiability class given in the Corollary \ref{convexinner} is the best possible result. 

\begin{theorem}
\label{differentiability}
Let $1\leq k \leq s \leq d$. Then there exists a polytope in ${\mathbb R}^d$ whose family ${\cal N}$ of inner normal vectors of facets has absolute rank $k$ and whose inner neighborhood volume function $V_P$ is of differentiability class ${\cal C}^{(s-1)}$, but not ${\cal C}^s$.
\end{theorem}

\begin{figure}[h]
\begin{center}
\includegraphics[scale=1.2]{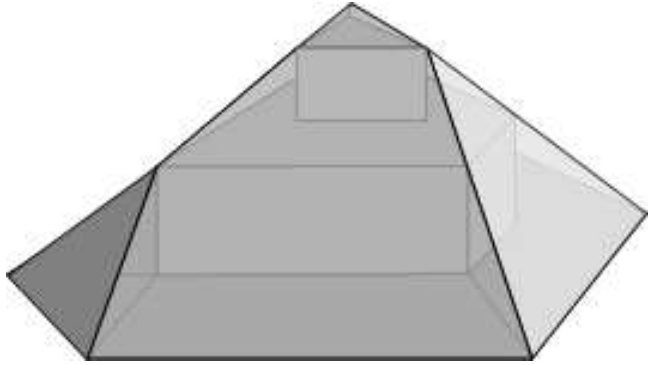}
\end{center}
\caption{ The roof of the polygon $ABCDE$ shown in Figure \ref{inner}} \label{roof}
\end{figure}

{\bf Proof} Consider the following construction, starting with be a polytope $P$ in ${\mathbb R}^d$ having $m$ facets, having inradius $g$ and absolute rank $k$ of the family ${\cal N}$ of normal vectors. We define the {\em roof of} $P$ to be the polytope $\Gamma(P)$ in  $R^{d+1}$ given by
$$\Gamma(P)=\{x'=(x,x_{d+1})| 0\leq x_{d+1}\leq g, x \in P(x_{d+1})\}.
$$
If, for all non-negative $r$, we canonically identify ${\mathbb R}^d$ with $\{(x,r)| x\in {\mathbb R}^d\}$, the $r$-interior of $P$ is congruent with the vertical slice $\{x_{d+1}=r\}$ in  $\Gamma(P)$, so $\Gamma(P)$ is the graph of the shrinking polytope $P_r$, as $r$ varies in $[0,+\infty)$ (see Figure \ref{roof}). 

Note that $\Gamma(P)$ is a polytope with $m+1$ facets and of absolute rank $k+1$. We identify $P$ with the facet $P \times \{0\}$ of $\Gamma(P)$, whose inner unit normal vector is the basis vector ${\mathbf e}_{d+1}$. Note that the outer dihedral angle between $P$ and any other facet of $\Gamma(P)$ equals $\frac{3\pi}{4}$. We conclude that, for all $r$, the translation of vector $r{\mathbf e}_{d+1}$ in $R^{d+1}$ maps $\Gamma\left(P((1+\sqrt{2})r)\right)$ onto $\Gamma(P)(r)$, thus the inradius of $\Gamma(P)$ is $\frac{g}{1+\sqrt{2}}$. As the volume of $\Gamma(P)$ is given by
$${\rm vol}_{d+1}(\Gamma(P))=\int_0^{\infty} {\rm vol}_{d}(P(\rho))\, d\rho
$$
we obtain that

$$\begin{array}{rl}
 {\rm vol}_{d+1}\left(\Gamma(P)(r)\right) & ={\rm vol}_{d+1}\left(\Gamma\left(P((1+\sqrt{2})r)\right)\right) \\ & \\
      & =\int_0^{\infty} {\rm vol}_{d}\left(P((1+\sqrt{2})r+\rho)\right)\, d\rho  \\ & \\
      & =\int_{(1+\sqrt{2})r}^{\infty} {\rm vol}_{d}(P(\rho))\, d\rho
\end{array}
$$
or, equivalently:
$$W'_{\Gamma(P)}(r)=-(1+\sqrt{2})W_P((1+\sqrt{2})r) \mbox{, for all $r\in[0,+\infty)$}.
$$

 Coming back to the proof of the theorem, let $d'=d-k+1$ and $m=s-k+1$ and consider a sequence $0< a_1= a_2 =  \dots =a_m < a_{m+1}\leq \dots a_{d'}$. Let $R$ be the rectangle $R_{a_1,a_2,\dots ,a_{d'} }$ in ${\mathbb R}^{d'}$ as in the example considered in the first section and let $P=\Gamma^{k-1}(R)=\Gamma(\Gamma(\dots(R)\dots))$. Then $P$ is a polytope in ${\mathbb R}^d$ of absolute rank $k$ and its volume function $W_P$ is of class ${\cal C}^{(s-1)}$, but not ${\cal C}^{s}$. qed.

{\em Aknowledgement}. We cordially  thank Ali Deniz and Yunus Ozdemir from Anadolu University for drawing the figures.


\begin{thebibliography}{99}
 \bibitem{Berger}
 M. Berger, {\it  Geometry}, Springer-Verlag, Berlin, 1987. 
 \bibitem{Morvan}
 J-M. Morvan, {\it  Generalized Curvatures}, Springer-Verlag, Berlin, 2008.
 \bibitem{Federer}
  H. Federer, {\it  Curvature measures}, Trans. Am. Math. Soc. vol. 93 (1959) pp. 418-491.
 \bibitem{Gray}
 A. Gray, {\it Tubes}, Birkhauser, 2004.
 \bibitem{Schneider}
 R. Schneider, {\it Convex Bodies: The Brunn-Minkowski Theory}, Cambridge Univ. Press, Cambridge, 1993. 
 \bibitem{LapidusF}
 M.L. Lapidus, M. van Frankenhuijsen, {\it Fractal Geometry, Complex Dimensions and Zeta Functions}, Springer-Verlag, Berlin, 2006. 
 \bibitem{LapidusP}
 M.L. Lapidus, E. P. J. Pearse, {\it Tube Formulas and Complex Dimensions of Self-Similar Tilings}, 
 Acta Mathematicae Applicandae, in press, 2009. 
 \bibitem{LapidusP1}
 M. L. Lapidus, E. P. J. Pearse, {\it Tube formulas for self-similar fractals}, in Analysis on Graphs 
and Its Applications (P. Exner, J. P. Keating, C. Bristol, P. Kuchment, T. Sunada, and A. Teplyaev, 
eds.), Proc. of Symposia in Pure Mathematics, vol. 77, Amer. Math. Soc., Providence, RI, 2008, 
pp. 211Ð230. arXiv:0711.0173.
  \bibitem{Deniz}
 A.Deniz, S.Kocak, Y. Ozdemir, A.E. Ureyen, {\it Tube formulas for self-similar fractals with non-Steiner-like generators}, available at arXiv: 0911.4966.
 \end{thebibliography}
 \end{document}